\input pictex.tex
\input amstex.tex
\NoBlackBoxes
\documentstyle{amsppt}
\magnification=1200
 \pagewidth{16.3truecm}
 \pageheight{24truecm}
 \nologo

\def\div{\operatorname{div}}
\def\mod{\operatorname{mod}}

\def\Pic{\operatorname{Pic}}
\def\Var{\operatorname{Var}}

\def\Var{\operatorname{Var}}

\def\deg{\operatorname{deg}}

\def\mot{\operatorname{mot}}
\def\card{\operatorname{card}}
\def\Hod{\operatorname{Hod}}
\def\top{\operatorname{top}}

\refstyle{A} \widestnumber\key{ACLM2}
 \topmatter
\title
Vanishing of principal value integrals on surfaces
\endtitle
\author
Willem Veys  
\endauthor
\address K.U.Leuven, Departement Wiskunde, Celestijnenlaan 200B,
         B--3001 Leuven, Belgium  \endaddress
\email wim.veys\@wis.kuleuven.be  \newline
 http://www.wis.kuleuven.be/algebra/veys.htm
\endemail
 \abstract
 Principal value integrals are associated to multi--valued
 rational differential forms with normal crossings support on
 a non--singular algebraic variety. We prove their vanishing on
 rational surfaces in the context of a conjecture of
 Denef--Jacobs. As an application we obtain a strong vanishing
 result for candidate poles of $p$--adic and motivic Igusa zeta
 functions.
\endabstract
 \endtopmatter

\document
\noindent {\bf Introduction}
 \bigskip
 \noindent {\bf 0.1.} Real
and $p$--adic principal value integrals were first introduced by
Langlands in the study of orbital integrals [Lan1,Lan2,LS1,LS2].
They are associated to multi--valued differential forms on real
and $p$--adic manifolds, respectively.

Let for instance $X$ be a non--singular projective algebraic
variety of dimension $n$ over $\Bbb Q_p$ (the field of $p$--adic
numbers). Denoting by $\Omega^n_X$ the vector space of {\it
rational} differential $n$--forms on $X$, take $\omega \in
(\Omega^n_X)^{\otimes d}$ defined over $\Bbb Q_p$; we then write
formally $\omega^{1/d}$ and consider it as a multi--valued
rational differential form on $X$.

We suppose that the support $|\div \omega|$ of $\div \omega$  has
normal crossings (over $\Bbb Q_p$) on $X$; say $D_i, i \in S$, are
its irreducible components. Let $\div \omega^{1/d} := \frac 1 d
\div \omega = \sum_{i \in S} (\alpha_i - 1) D_i$, where then the
$\alpha_i \in \frac 1 d \Bbb Z$. If $\omega^{1/d}$ has {\it no
logarithmic poles}, i.e. if all $\alpha_i \ne 0$, the {\it
principal value integral} $PV \int_{X(\Bbb Q_p)} |
\omega^{1/d}|_p$ of $\omega^{1/d}$ on $X(\Bbb Q_p)$ is defined as
follows.  Cover $X(\Bbb Q_p)$ by (disjoint) small enough open
balls $B$ on which there exist local coordinates $x_1, \dots ,
x_n$ such that all $D_i$ are coordinate hyperplanes. Consider for
each $B$ the {\sl converging} integral $\int_B |x_1 x_2 \cdots x_n
|^s_p | \omega^{1/d}|_p$ for $s \in \Bbb C$ with ${\Cal R}(s) >\!>
0$, take its meromorphic continuation to $\Bbb C$ and evaluate
this in $s = 0$; then add all these contributions. One can check
that the result is independent of all choices.

In the real setting we proceed similarly but then we also need a
partition of unity, and we have to assume that $\omega^{1/d}$ has
no integral poles, i.e. the $\alpha_i \notin \Bbb Z_{\leq 0}$.
Here the independency result is somewhat more complicated; it was
verified in detail in [Ja1].

\bigskip
\noindent {\bf 0.2.} Denef and Jacobs proved a vanishing result
for real principal value integrals, and conjectured a similar
statement in the $p$--adic case. In both cases let ${\Cal
L}(\omega^{1/d})$ be the locally constant sheaf of $\Bbb
C$--vector spaces on $X \setminus | \div \omega |$ associated to
$\omega^{1/d}$. It has rank 1, a non--zero section on a connected
open  being an analytic branch of $\omega^{1/d}$ multiplied with a
complex number. (In the $p$--adic case we choose an embedding of
$\Bbb Q_p$ into $\Bbb C$.)

\bigskip
\proclaim{Theorem {\rm [DJ, 1.1.4],[Ja1]}} Let $X$ be a
non--singular projective algebraic variety of \linebreak (complex)
dimension $n$, defined over $\Bbb R$. If $H^n(X(\Bbb C) \setminus
| \div \omega |, {\Cal L}(\omega^{1/d})) = 0$, then \linebreak $PV
\int_{X(\Bbb R)} | \omega^{1/d} | = 0$. \endproclaim

\bigskip
\proclaim{Conjecture {\rm[DJ, 1.2.2]}} Let $X$ be a non-singular
projective algebraic variety, defined over $\Bbb Q_p$. If $H^i
(X(\Bbb C) \setminus | \div \omega |, {\Cal L} (\omega^{1/d})) =
0$ for all $i \geq 0$, then $PV \int_{X(\Bbb Q_p)} | \omega^{1/d}
|_p = 0$.
\endproclaim

\medskip
 \noindent (The authors are cautious and mention that
perhaps one has to suppose also some good reduction $\mod p$ and
that all $\alpha_i \not\in \Bbb Z$.)

\bigskip
\noindent {\bf 0.3.} In [Ve7] we `upgraded' $p$--adic principal
value integrals to motivic ones, in the same spirit as how motivic
integration and motivic zeta functions were inspired by (usual)
$p$--adic integration and $p$--adic Igusa zeta functions. See
[DL2,DL3] or the surveys [DL4,Loo,Ve6] for these notions.

More precisely, let $X$ be a non--singular algebraic variety (say
over $\Bbb C$) of dimension $n$ and $\omega^{1/d}$ a multi--valued
differential form on $X$. Let as above $\div(\omega^{1/d}) =
\sum_{i \in S} (\alpha_i - 1)D_i$ be a normal crossings divisor.
We denote also $D^\circ_I := (\bigcap_{i \in I} D_i) \setminus
(\bigcup_{\ell \not \in I} D_\ell)$ for $I \subset S$; so $X =
\coprod_{I \subset S} D^\circ_I$. If $\omega^{1/d}$ has no
logarithmic poles, then the {\it motivic principal value integral}
of $\omega^{1/d}$ on $X$ is
$$PV \int_X \omega^{1/d} = L^{-n} \sum_{I \subset S} [D^\circ_I]
\prod_{i \in I} \frac{L-1}{L^{\alpha_i} - 1}.$$
 Here $[\cdot]$
denotes the class of a variety in the Grothendieck ring of
algebraic varieties, and $L := [ \Bbb A^1]$, see (1.1). We refer
to [Ve7, \S 2] for a motivation for this expression. Note for
instance that this is just the formula for the {\it converging}
motivic integral associated to the $\Bbb Q$--divisor
$\div(\omega^{1/d})$ if all $\alpha_i > 0$.

It is natural to `upgrade' also the conjecture in (0.2) to this
setting (maybe also assuming that all $\alpha_i \not\in \Bbb Z$) :
{\sl if $H^i(X(\Bbb C) \setminus | \div \omega | , {\Cal
L}(\omega^{1/d})) = 0$ for all $i \geq 0$, then $PV \int_X
\omega^{1/d} = 0$.}

\bigskip
\noindent {\bf 0.4.} In this paper we attack (the motivic version
of) the conjecture for surfaces, more precisely {\it rational}
surfaces.

We note that, in dimension 2, this is the crucial class to study :
there is no `classification' for configurations $| \div \omega |
\subset X$ with all $H^i (X \setminus | \div \omega |, {\Cal
L}(\omega^{1/d})) = 0$ on a rational surface $X$. For
non--rational surfaces it is conceivable that the conjecture can
be approached through the classification of such configurations
from [GP] and [Ve3]. We plan to report on this later. (Also for
the applications that we  will prove here, the class of rational
surfaces is the crucial one, see (0.7).)

Our main theorem is as follows.

\bigskip
\proclaim{Theorem} Let $X$ be a non--singular projective rational
surface and $\omega^{1/d}$ a multi--valued differential form on
$X$ without logarithmic poles (in particular $|\div \omega|$ has
normal crossings).

(1) Suppose that $B := | \div \omega |$ is {\rm connected}. If
$\chi(X \setminus B) \leq 0$, then $PV \int_X \omega^{1/d} = 0$.

(2) More generally, let $B$ be any {\rm connected} normal
crossings divisor satisfying $B \supset | \div \omega |$. If $\chi
(X \setminus B) \leq 0$, then $PV \int_X \omega^{1/d} = 0$.
\endproclaim

\bigskip
\noindent Statement (1) is somewhat weaker than the conjecture
because of the connectivity condition. On the other hand it is
clearly stronger since we assume only that $\chi (X \setminus B)
\leq 0$, instead of the vanishing of all $H^i (X \setminus B,
{\Cal L}(\omega^{1/d}))$. (And we do not need the extra assumption
$\alpha_i \not\in \Bbb Z$.) The generalization (2) is important in
view of the applications on motivic zeta functions and is natural
in that context, see (0.8).

The important ingredients in our proof are the structure theorem
of [Ve3] for such configurations $B \subset X$ with $\chi (X
\setminus B) \leq 0$, and the notion of a more general principal
value integral when $\omega^{1/d}$ is allowed to have `some'
logarithmic poles, see (1.4).

\bigskip
\noindent {\bf 0.5.} Real and $p$--adic principal value integrals
appear as coefficients of asymptotic expansions of oscillating
integrals and fibre integrals, and as residues of poles of
distributions $|f|^\lambda$ and $p$--adic Igusa zeta functions,
respectively. See [Ja2, \S 1] for an overview and
[AVG,De2,Ig1,Ig2,Ja3,Lae] for more details.

We have in particular that the cancelation of some given candidate
pole of a $p$--adic Igusa zeta function is essentially equivalent
to the vanishing of an associated $p$--adic principal value
integral. Here we are interested in the analogous phenomenon for
poles of the motivic zeta function versus associated motivic
principal value integrals.

\bigskip
\noindent{\bf 0.6.} Denef and Loeser [DL2] associated to a
non--constant regular function $f$ on a non--singular algebraic
variety $M$ of dimension $n+1$ its motivic zeta function ${\Cal
Z}_{\mot}(f;T)$; here $T$ is a variable. They obtained the
following formula for it in terms of an embedded resolution $h : Y
\rightarrow M$ of the hypersurface $\{ f = 0 \}$. Let $E_j, j \in
K$, be the irreducible components of $h^{-1} \{ f = 0 \}$ and
$N_j$ and $\nu_j - 1$ the multiplicity of $E_j$ in the divisors
$\div(f \circ h)$ and $\div(h^\ast dx)$, respectively, where $dx$
is a local generator of the sheaf of $(n+1)$--forms on $M$. Denote
$E^\circ_J := (\cap_{j \in J} E_j) \setminus (\cup_{\ell \not\in
J} E_\ell)$ for $J \subset K$. Then
$${\Cal Z}_{\mot} (f;T) = L^{-(n+1)} \sum_{J \subset K} [E^\circ_J]
\prod_{j \in J} \frac{(L-1)T^{N_j}}{L^{\nu_j} - T^{N_j}}\, .$$
 Fix a `candidate pole' $L^{\nu_j/N_j}$ of
${\Cal Z}_{\mot} (f;T)$; see (3.4) for more explanations. In the
generic situation that $\nu_j/N_j \ne \nu_i/N_i$ for all $i \ne
j$, the cancelation of the candidate pole (of order 1)
$L^{\nu_j/N_j}$ is equivalent to the vanishing of `its residue'
$$L^{-n} \sum_{I \subset S_j} [D^\circ_I] \prod_{i \in I}
\frac{L-1}{L^{\alpha_i}-1},$$ where $S_j = \{ i \in K \setminus \{
j \} \mid E_i \text{ intersects } E_j \}$, $D_i := E_j \cap E_i$
and $\alpha_i = \nu_i - (\nu_j/N_j) N_i \neq 0$ for $i \in S_j$,
and $D^\circ_I$ is as usual. This looks like a motivic principal
value integral on $E_j$, and indeed it is equal to $PV \int_{E_j}
\omega^{1/d}$, where $\omega^{1/d}$  is some Poincar\'e residue,
see (3.5).

Suppose now that $E_j$ is mapped to a point by $h$. Denote $B_j :=
\cup_{i \in S_j} D_i$. The famous Monodromy Conjecture predicts
more or less that (generically), if $(-1)^n \chi (E_j \setminus
B_j) \leq 0$, then $L^{\nu_j/N_j}$ is no pole of $Z_{\mot} (f;T)$,
see (3.5).

\bigskip
\noindent {\bf 0.7.} When $n=2$, we now can prove this expected
cancelation of candidate poles as a consequence of our Main
Theorem. More precisely, we may suppose that $h$ is a composition
of blowing--ups as in [Hi]. Then the exceptional surface $E_j$ is
created during the resolution process $h$ either as a projective
plane by blowing up a point, or as a ruled surface by blowing up a
non--singular curve. Then, with $B_j$ as above, we can show as
corollary of (0.4):

\bigskip
\proclaim{Theorem} Let $L^{\nu_j/N_j}$ be a candidate pole of
order 1 for ${\Cal Z}_{\mot} (f;T)$, as described above. Suppose
that $\chi(E_j \setminus B_j) \leq 0$.

(1) If $E_j$ is created by blowing up a point, then we have {\rm
always} that $L^{\nu_j/N_j}$ is no pole.

(2) If $E_j$ is created by blowing up a rational curve, then
$L^{\nu_j/N_j}$ is no pole whenever $B_j$ is connected.
\endproclaim

\medskip \noindent In fact we obtain a somewhat stronger statement;
see (0.8(i)) below. But first we want to comment on this new
result.

Concerning (1), the best general result up to now [Ve2] was the
analogous cancelation in the context of $p$--adic Igusa zeta
functions when the centre of the blowing up is a point of
multiplicity at most 4 on the strict transform of $\{ f = 0 \}$.
This was achieved using a lengthy classification of all possible
configurations with $\chi(E_j \setminus B_j) \leq 0$ under this
multiplicity restriction.

Concerning (2), it is important to note that non--connected
intersection configurations on $E_j$ are very rare, see (3.7). And
moreover, the case of {\it rational} centres is the crucial one.
When $E_j$ is created by blowing up a curve of genus $g \geq 2$ we
already proved the expected cancelation in [Ve2], and we can now
also handle the case $g=1$ completely by combining [Ve2] and [Ve3]
with recent work of Rodrigues [Ro2].

\bigskip
\noindent {\bf 0.8.} {\sl Remarks.} (i) In fact we assume only
that $\nu_j/N_j \ne \nu_i/N_i$ for the $i \in S_j$, and we show
that `the contribution of $E_j$ to the residue of $L^{\nu_j/N_j}$'
vanishes when expected.

(ii) Similar vanishing results hold in the context of much more
general motivic zeta functions, for instance those of [Ve4]
associated to an effective $\Bbb Q$--Cartier divisor on any $\Bbb
Q$--Gorenstein threefold (instead of just a hypersurface on $\Bbb
A^3$). Also, (0.7) specializes to the context of Hodge,
topological and $p$--adic zeta functions.

(iii) Note that in (0.6) it is possible that some $\alpha_i = 1$,
and thus that $B_j (:= \cup_{i \in S_j} D_i) \supsetneqq |\div
\omega^{1/d}|$. So the more general setting of part (2) of our
main theorem pops up naturally in this context of poles of zeta
functions.

\bigskip
\noindent {\bf 0.9.} We will work over the base field $\Bbb C$ of
complex numbers. In \S1 we introduce the more general principal
value integrals on surfaces, allowing \lq some\rq\ logarithmic
poles, which we need for the proof of the main theorem in \S2.
Then in \S3 we obtain as a corollary the cancelation of candidate
poles for zeta functions.

\bigskip
\noindent {\sl Acknowledgement.} We would like to remark that
[ACLM1, \S2] contributed to our inspiration for the proof of the
main theorem. And we thank Bart Rodrigues for his useful remarks
and suggestions.

\bigskip
\bigskip
\noindent {\bf 1. Generalized principal value integrals on
surfaces}

\bigskip
\noindent {\bf 1.1.}   We first recall briefly the notion of
Grothendieck ring of varieties and related constructions.

(i) The Grothendieck ring $K_0 (\Var)$ of complex algebraic
varieties  is the free abelian group generated by the symbols
$[V]$, where $V$ is a variety, subject to the relations $[V] =
[V^\prime]$ if $V$ is isomorphic to $V^\prime$, and $[V] = [V
\setminus W] + [W]$ if $W$ is closed in $V$. Its ring structure is
given by $[V] \cdot [W] := [V \times W]$. (This ring is quite
mysterious; see [Po] for the recent proof that it is not a
domain.) Usually, one abbreviates $L := [\Bbb A^1]$.

For the sequel we need to extend $K_0(\Var)$ with fractional
powers of $L$ and to localize. Fix $d \in \Bbb Z_{> 0}$; we
consider
$$K_0 (\Var)[L^{-1/d}] := {K_0(\Var)[T] \over (LT^d - 1)}$$
(where $L^{-1/d} := \bar T$). We then localize this ring with
respect to the elements $L^{i/d} - 1$,
  $i \in \Bbb Z
\setminus \{ 0 \}$. What we really need is the subring of this
localization generated by $K_0 (\Var)$, $L^{-1}$ and the elements
$(L-1)/(L^{i/d} - 1), i \in \Bbb Z \setminus \{ 0 \}$; we denote
this subring by ${\Cal R}_d$. (We do not know whether or not
${\Cal R}_d$ has zero divisors.)

(ii) For a variety $V$, we denote by $h^{p,q} (H^i_c (V, \Bbb C))$
the rank of the $(p,q)$--Hodge component in the mixed Hodge
structure of the $i$th cohomology group with compact support of
$V$. The {\it Hodge polynomial} of $V$ is
$$H(V) = H(V;u,v) := \sum_{p,q} \big(\sum_{i \geq 0} (-1)^i h^{p,q}
(H^i_c (V,\Bbb C))\big) u^pv^q \in \Bbb Z [u,v].$$
 Precisely by
the defining relations of $K_0 (\Var)$, there is a well--defined
ring homomorphism $H : K_0 (\Var) \rightarrow \Bbb Z[u,v]$,
determined by $[V] \mapsto H(V)$. It induces a ring homomorphism
$H$ from ${\Cal R}_d$ to the `rational functions in $u,v$ with
fractional powers'.

(iii) The topological Euler characteristic $\chi(V)$ of a variety
$V$ satisfies $\chi(V)=H(V;1,1)$ and we obtain a ring homomorphism
$\chi : K_0 (\Var) \rightarrow \Bbb Z$, determined by $[V] \mapsto
\chi(V)$. Since $\chi(L)=1$, it induces a ring homomorphism $\chi:
{\Cal R}_d \to \Bbb Q$ by declaring $\chi((L-1)/(L^{i/d} -
1))=d/i$ (see e.g. [DL2,Ve4] for similar constructions).

\bigskip
\noindent {\bf 1.2.} On surfaces we want to extend the notion of
principal value integral of (0.3) in two ways. First we will allow
the differential form $\omega^{1/d}$ to have `some' logarithmic
poles. Another somewhat technical generalization consists in
considering a normal crossings divisor whose support {\sl
contains} $| \div \omega |$. In view of the application on
candidate poles of zeta functions, this is natural; see (0.8(iii))
or (3.5). More precisely we fix the following setting.

\bigskip
\noindent {\bf 1.3.} Let $X$ be a non--singular projective
surface, $\omega^{1/d}$ a multi--valued differential form on $X$,
and $D = \cup_{i \in T} C_i$ a normal crossings divisor on $X$
satisfying $D \supset | \div \omega |$. We write formally $\div
\omega^{1/d} = \sum_{i \in T} (\alpha_i - 1) C_i$, where thus
$\alpha_i = 1$ if $C_i \not\subset | \div \omega |$. We put the
following restriction on the possible logarithmic poles of
$\omega^{1/d}$.

\noindent{\sl If for $i \in T$ we have $\alpha_i = 0$, then

(1) $C_i$ is rational,

(2) no $C_\ell$ that intersects $C_i$ has $\alpha_\ell  = 0$,

(3) in all but at most two intersection points of $C_i$ with other
$C_\ell$ the intersecting curve $C_\ell$ has $\alpha_\ell = 1$.}

\noindent We call such a pair $(D, \omega^{1/d})$ {\sl allowed}.

\pagebreak

\bigskip
When $\alpha_i = 0$ the adjunction formula $(K_X + C_i) \cdot C_i
= K_{C_i}$ and the restrictions above easily yield that either two
curves $C_{i_1}$ and $C_{i_2}$ with $\alpha_{i_1} \ne 1 \ne
\alpha_{i_2}$ intersect $C_i$ and then $\alpha_{i_1} +
\alpha_{i_2} = 0$, or only one curve $C_{i_1}$ with $\alpha_{i_1}
\ne 1$ intersects $C_i$ and then $\alpha_{i_1} = -1$. (Here and
further on, when applying the adjunction formula we always
consider $\div \omega^{1/d}$ as a representative of $K_X$.)

\bigskip
\noindent {\bf 1.4. Definition.} Let $X$ be a non--singular
projective surface and $(D,\omega^{1/d})$ an allowed pair on $X$.
We write $\div \omega^{1/d} = \sum_{i \in T} (\alpha_i - 1)C_i$ as
in (1.3) and $C^\circ_I := (\bigcap_{i \in I} C_i) \setminus
(\bigcup_{\ell \not\in I} C_\ell)$ for $I \subset T$ as before.
Furthermore for $i \in T$     
 we denote by
$C_i \cdot C_i$ the self--intersection number of $C_i$ and by
$C_j, j \in T_i (\subset T)$, the curves that intersect $C_i$. To
$(D,\omega^{1/d})$ we associate the invariant
$${\Cal E}_X (D,\omega^{1/d}) := \sum \Sb I \subset T \\ \forall i
\in I : \alpha_i \ne 0 \endSb [C^\circ_I] \prod_{i \in I}
\frac{L-1}{L^{\alpha_i} - 1} + \sum \Sb i \in T \\ \alpha_i = 0
\endSb (-C_i \cdot C_i) \prod_{j \in T_i} \frac{L-1}{L^{\alpha_j}
- 1},
$$
living in ${\Cal R}_d$. In the last sum the expression $(-C_i
\cdot C_i) \prod_{j \in T_i} \frac{L-1}{L^{\alpha_j}-1}$ is the
easiest to write down `uniformly', but, since at most two of the
$\alpha_j$ are different from 1, this expression boils down to the
following.

(1) If $C_i$ intersects two curves $C_{i_1}$ and $C_{i_2}$ with
$\alpha_{i_1} \ne 1 \ne \alpha_{i_2}$, we get
$$(-C_i \cdot C_i)
\frac{(L-1)^2}{(L^{\alpha_{i_1}}-1)(L^{\alpha_{i_2}} - 1)} = (C_i
\cdot C_i) \frac{(L-1)^2 L^{\alpha_{i_1}}}{(L^{\alpha_{i_1}} -
1)^2} = (C_i \cdot C_i) \frac{(L-1)^2
L^{\alpha_{i_2}}}{(L^{\alpha_{i_2}}-1)^2}\, .$$

(2) If $C_i$ intersects only one curve $C_{i_1}$ with
$\alpha_{i_1} \ne 1$, we get
$$(-C_i \cdot C_i) \frac{L-1}{L^{\alpha_{i_1}}-1} = (C_i \cdot
C_i)L\, .$$

\noindent {\sl Note.} (i) In fact the $C_i$ with $\alpha_i = 1$,
i.e. those $C_i$ not belonging to $|\div \omega |$, play no role
in the definition of ${\Cal E}_X(D,\omega^{1/d})$: we could as
well consider instead of $T$ only $\{ i \in T \mid C_i \subset |
\div \omega | \}$. So this invariant is really an invariant of
$\omega^{1/d}$ only. However, for the sequel it is useful to
introduce it as above.

(ii) As a motivation for the expression for the contribution of
$C_i$ with $\alpha_i = 0$: it is a kind of limit of the `total
contribution of $C_i$' in the formula of (0.3) for $PV \int_X
\omega^{1/d}$ if $\alpha_i \ne 0$, see e.g. [Ve5, 3.3].

\vskip 1truecm
 \centerline{
\beginpicture
\setcoordinatesystem units <.5truecm,.5truecm>
 \putrule from 1 2
to 11 2
 \putrule from 2.5 1 to 2.5 5
  \putrule from 4 1 to 4 5
 \setdashes
  \putrule from 5.5 1 to 5.5 5
  \putrule from 7.5 1 to 7.5 5
  \putrule from 9.5 1 to 9.5 5
 \setsolid
\put {$\dots$} at 6.5 4.5
 \put {$C_i$} at 1.1 2.5
 \put {$C_{i_1}$} at 3.1 5.3
 \put {$C_{i_2}$} at 4.6 5.3
 \put {$C_j$} at 10.2 5
 \put {$\bullet$} at 9.5 2
 \put {$P$} at 10 2.5
\put {$\longleftarrow$} at 14 2.8
 \put {$h$} at 14 3.5
\setcoordinatesystem units <.5truecm,.5truecm> point at -16 0
\putrule from 1 2 to 11 2
 \putrule from 2.5 1 to 2.5 5
  \putrule from 4 1 to 4 5
  \setlinear \plot 9 1  10.5 4 /
 \setdashes
  \putrule from 5.5 1 to 5.5 5
  \putrule from 7.5 1 to 7.5 5
  \plot  8.8 5    10.9 2.9 /
   \setsolid
\put {$\dots$} at 6.5 4.5
 \put {$C_i$} at 1.1 2.5
 \put {$C_{i_1}$} at 3.1 5.3
 \put {$C_{i_2}$} at 4.6 5.3
 \put {$C_j$} at 9.7 5.1
 \put {$\bullet$} at 9.5 2
 \put {$C$} at 9.5 2.9
\endpicture
 }
 \vskip 1truecm
  \centerline{\smc Figure 1}
 \vskip 1truecm

\proclaim{1.5. Lemma} Let $X$ be a non--singular projective
surface, $P \in X$ and $h : \tilde X \rightarrow X$ the
blowing--up of $X$ with centre $P$. Let $(D,\omega^{1/d})$ be
allowed on $X$, writing as usual $\div \omega^{1/d} = \sum_{\ell
\in T} (\alpha_\ell - 1)C_\ell$. Then

(1) $(h^{-1}D, h^\ast \omega^{1/d})$ is allowed on $\tilde X$,

(2) ${\Cal E}_{\tilde X} (h^{-1} D, h^\ast \omega^{1/d}) = {\Cal
E}_X (D,\omega^{1/d})$ except when (on $X$) there exists a curve
$C_i$ with $\alpha_i = 0$ and curves $C_{i_1}$ and $C_{i_2}$,
intersecting $C_i$, with $\alpha_{i_1}+\alpha_{i_2}=0$ and $\{
\alpha_{i_1}, \alpha_{i_2} \} \ne \{ -1,1 \}$, such that $P \in
C_i$, $P \not\in C_{i_1}$ and $P \not \in C_{i_2}$. (And in this
exceptional case we do have inequality.)
\endproclaim

\bigskip
\demo{Proof} The necessary (easy) computations are essentially in
[ACLM1] and [Ve5, 3.5]. We just illustrate the exceptional case,
see Figure 1. Note that in this case, since $(D, \omega^{1/d})$ is
allowed, if there is another curve $C_j$ passing through $P$, it
must satisfy $\alpha_j = 1$. And by Note (1.4(i)) we may as well
assume that only $C_{i_1}$ and $C_{i_2}$ intersect $C_i$.

Let $C$ denote the exceptional curve of $h$. Then $h^{-1}D$
consists of the union of $C$ and the strict transforms of the
$C_\ell$. And since $\alpha_i = 0$ we have that $C$ does not
appear in the divisor of $h^\ast \omega^{1/d}$. So we write
formally $\div(h^\ast \omega^{1/d}) = \sum_{\ell \in T}
(\alpha_\ell - 1) C_\ell + (\alpha - 1)C$ with $\alpha = 1$. We
have to compare the contributions of $C_i$ to ${\Cal E}_X (D,
\omega^{1/d})$ and of $C_i \cup C$ to ${\Cal E}_{\tilde X}
(h^{-1}D, h^\ast \omega^{1/d})$. These are
$$-(C_i \cdot C_i) \frac{(L-1)^2}{(L^{\alpha_{i_1}} -
1)(L^{\alpha_{i_2}} - 1)}$$ and
$$-(C_i \cdot C_i - 1) \frac{(L-1)^2}{(L^{\alpha_{i_1}} -
1)(L^{\alpha_{i_2}} - 1)} + L ,$$ respectively. Their difference
$\frac{(L-1)^2}{(L^{\alpha_{i_1}} - 1)(L^{\alpha_{i_2}} - 1)} + L$
is nonzero (in ${\Cal R}_d$); one can ensure the non--nullity for
instance using $H$ or $\chi$. (On the other hand, when $\{
\alpha_{i_1}, \alpha_{i_2} \} = \{ -1,1 \}$, this difference would
be $-L + L = 0$.) \qed
\enddemo

\bigskip
\bigskip
\noindent {\bf 2. Rational surfaces}

\bigskip
\noindent {\bf 2.1.} We first summarize the structure theorem of
[Ve3] and some of its refinements, which will be the starting
point of the proof of our main theorem. Remember that a
non--singular rational curve with self--intersection $-1$ is
called a $(-1)$--curve.

\bigskip
\noindent \proclaim{Structure Theorem} Let $X$ be a non--singular
projective rational surface. Let $B$ be a {\rm connected} normal
crossings curve on $X$ with $\chi(X \setminus B) \leq 0$. Assume
that $X$ does not contain any $(-1)$--curve disjoint from $B$.

By [GP, Theorem 3] there is a dominant morphism $\varphi : X
\setminus B \rightarrow \Bbb P^1$; let $h : \tilde X \rightarrow
X$ be the {\rm minimal} morphism that resolves the indeterminacies
of $\varphi$, considered as rational map from $X$ to $\Bbb P^1$.

(1) Then there exists a connected curve $B^\prime \supset B$ with
$\chi (X \setminus B^\prime) \leq \chi(X \setminus B) \leq 0$,
such that the morphism $\varphi \circ h$ decomposes as
$$\tilde X \overset g \to \longrightarrow \Sigma \overset \pi \to
\longrightarrow \Bbb P^1,$$ where $g$ is a composition of
blowing--downs with exceptional curve in $h^{-1}B^\prime$, and
$\pi : \Sigma \rightarrow \Bbb P^1$ is a ruled surface; see
Diagram 1. Moreover, $h^{-1}B^\prime$ has normal crossings in
$\tilde X$.

(2) We can require the configuration $g(h^{-1}B^\prime) \subset
\Sigma$ to be one of the configurations in Figure 2.

\vskip 1truecm
 \centerline{
\beginpicture
 \setcoordinatesystem units <.5truecm,.5truecm> \putrectangle
corners at 0 0 and 10 6 \putrule from 1 3 to 9 3 \putrule from 2.5
1 to 2.5 5 \putrule from 4 1 to 4 5 \setdashes \putrule from 5.5 1
to 5.5 5 \putrule from 8 1 to 8 5 \setsolid
 \put {$\dots$} at 6.75
4.5 \put {$C_1$} at 1.2 2.3
 \put {(a)} at 5 -1.3
 \setcoordinatesystem units <.5truecm,.5truecm> point at -11 0
\putrectangle corners at 0 0 and 10 6 \putrule from 1 2 to 9 2
\putrule from 1 4 to 9 4 \putrule from 3 1 to 3 5 \setdashes
\putrule from 5 1 to 5 5 \putrule from 8 1 to 8 5 \setsolid
 \put
{$\dots$} at 6.5 3 \put {$C_1$} at 1.2 3.3 \put {$C_2$} at 1.2 1.3
\put {(b)} at 5 -1.3
 \setcoordinatesystem units
<.5truecm,.5truecm> point at -22 0
 \putrectangle corners at 0 0
and 10 6 \putrule from 1 1 to 1 5 \putrule from 9 1 to 9 5
\setdashes \putrule from 3 1 to 3 5 \putrule from 5 1 to 5 5
\setsolid
 \ellipticalarc axes ratio 4:1  360 degrees from 1 3
center at 5 3
 \put {$\dots$} at 4 4.75 \put {$C$} at 7 4.5
\put {(c)} at 5 -1.3
\endpicture
 }
 \vskip 1truecm
  \centerline{\smc Figure 2}
 \vskip 1truecm

\noindent Here $C_1$ and $C_2$ are sections of $\pi$, $C$ is a
non--singular curve for which $\pi|_C : C \rightarrow \Bbb P^1$
has degree 2 (a `bisection'), and the other curves are fibres of
$\pi$. The minimal number of fibres in (a) and (b) is 2 and 1,
respectively; in (c) there must pass a fibre through each
ramification point of $\pi|_C$, and we can have any number of
other fibres. Note that in (c) the bisection can be non--rational
(and then has more than two ramification points).

(3) Irreducible curves ${\Cal C} \subset h^{-1}B^\prime$, which
are not components of $h^{-1}B$, occur only in fibres ($=$
exceptional components) of $g$. Moreover, any fibre of $g$
contains {\rm at most one} such curve ${\Cal C}$  and
$$\cases \card ({\Cal C} \cap h^{-1}B) = 1 & \text{if}\quad  \chi(\Sigma
\setminus g(h^{-1} B^\prime)) < 0 \\
\card ({\Cal C} \cap h^{-1}B) = 2 & \text{if} \quad \chi(\Sigma
\setminus g(h^{-1}B^\prime)) = 0. \endcases$$
 Note that always $\chi(\Sigma \setminus g(h^{-1}
B^\prime)) = 0$ in cases (b) and (c). \endproclaim

\bigskip
\noindent \demo {Proof} Combine essentially (3.3), (3.5) and (4.3)
in [Ve3]. \qed \enddemo

\vskip 1truecm
 \centerline{
\beginpicture
\setcoordinatesystem units <.5truecm,.5truecm>
 \put{$\searrow$} at
1 1 \put{$\swarrow$} at 1 -1 \put{$\searrow$} at -1 -1
\put{$\swarrow$} at -1 1 \put{$\Bigg\downarrow$} at 0 0
 \put{$g$}
at -1.3 1.3
 \put{$h$} at 1.3 1.3
 \put{$\tilde \varphi$} at -.4 0
 \put{$\varphi$} at 1.3 -1.3
 \put{$\pi$} at -1.3 -1.3
\put{$\Sigma$} at -2 0
 \put{$\tilde X$} at  0 2
 \put{$X$} at 2 0
 \put{${\bold P}^1$} at  0 -2
\endpicture
 }
 \vskip 1truecm
 \centerline{\smc Diagram 1}
 \vskip 1truecm

\pagebreak

\noindent {\bf 2.2.} We will denote by $p : {\tilde \Sigma}
\rightarrow \Sigma$ the minimal embedded resolution of the
configuration $D := g(h^{-1} B^\prime) \subset \Sigma$ in case (c)
of the structure theorem. When $\omega^{1/d}$ is a multi--valued
differential form on $\Sigma$ with $D \supset | \div \omega |$, we
will slightly abuse the terminology of (1.3) and say that
$(D,\omega^{1/d})$ is allowed on $\Sigma$ if $(p^{-1}D,p^\ast
\omega^{1/d})$ is allowed on $\tilde \Sigma$. In that case we also
put ${\Cal E}_{\Sigma} (D, \omega^{1/d}) := {\Cal E}_{\tilde
\Sigma} (p^{-1} D, p^\ast \omega^{1/d})$.

\bigskip
\noindent \proclaim {Lemma} Let $D := g(h^{-1} B^\prime) \subset
\Sigma$ in the Structure Theorem and let $\omega^{1/d}$ be a
multi--valued differential form on $\Sigma$. Assume in case (c)
that $C$ is rational (this is equivalent to $\pi|_C$ having
exactly two ramification points). If $(D, \omega^{1/d})$ is
allowed on $\Sigma$, then ${\Cal E}_{\Sigma} (D,\omega^{1/d})=0$.
\endproclaim

\medskip
 \demo{Proof} We write as usual $D = \bigcup_{i \in T} C_i$ and
$\div \omega^{1/d} = \sum_{i \in T} (\alpha_i - 1) C_i$. Applying
the adjunction formula to a generic fibre of $\pi$ yields in cases
(a), (b) and (c) that $\alpha_1 = -1$, $\alpha_1 + \alpha_2 = 0$
and $\alpha = 0$, respectively. We first treat the cases (a) and
(b).

 If no $\alpha_i = 0$ this is well known and easily
verified; see e.g. [Ve2] for a similar computation. The point is
that the contribution of {\it any} fibre of $\pi$ to ${\Cal E}_X
(D, \omega^{1/d})$ is zero. Now if $\alpha_i = 0$ for some fibre
$C_i$ of $\pi$, the contribution of $C_i$ is still zero, simply
because its self--intersection $C_i \cdot C_i = 0$. This finishes
already case (a).

In case (b) we are left with the following possibility: $\alpha_1
= \alpha_2 = 0$ and (omitting the possible fibres $C_\ell$ with
$\alpha_\ell=1$) there is either only one fibre $C_i \subset D$
with then necessarily $\alpha_i = -1$, or there are two fibres
$C_i$ and $C^\prime_i$ in $D$ with $\alpha_i + \alpha_i^\prime =
0$ and $\alpha_i \ne 0 \ne \alpha^\prime_i$. We compute ${\Cal
E}_X(D,\omega^{1/d})$ in this last case:
$$\split {\Cal E}_X(D,\omega^{1/d}) = & (L-1)^2 + (L-1)
\frac{L-1}{L^{\alpha_i}-1} + (L-1) \frac{L-1}{L^{\alpha^\prime_i}
- 1} \\
& + (- C_1 \cdot C_1)
\frac{(L-1)^2}{(L^{\alpha_i}-1)(L^{\alpha^\prime_i}-1)} + (-C_2
\cdot C_2)
\frac{(L-1)^2}{(L^{\alpha_i}-1)(L^{\alpha^\prime_i}-1)}.
\endsplit
$$
 Since $C_1 \cdot C_1 = -C_2 \cdot C_2$ (see e.g. [Ha, Theorem
 V.2.17]) the last two terms cancel and, since $\alpha_i +
\alpha^\prime_i = 0$, we obtain
$${\Cal E}_X(D,\omega^{1/d}) = (L-1)^2 \frac{L^{\alpha_i +
\alpha_i^\prime} - 1}{(L^{\alpha_i}-1)(L^{\alpha^\prime_i}-1)} =
0.$$

For case (c) we have to consider ${\Cal E}_{{\tilde \Sigma}}
(p^{-1} D, p^\ast \omega^{1/d})$.  We denote the two curves in
$p^{-1}D$ which intersect $C(\subset {\tilde \Sigma})$ in $P_1$
and $P_2$ by $C_1$ and $C_2$, respectively, see Figure 3.

\vskip 1truecm
 \centerline{
\beginpicture
\setcoordinatesystem units <.5truecm,.5truecm>
 \putrectangle
corners at -.5 0 and 10.5 6
 \putrule from 1 1 to 1 5
 \putrule from 9 1 to 9 5
 \setdashes
 \putrule from 4 1 to 4 5
 \setsolid
\ellipticalarc axes ratio 4:1  360 degrees from 1 3 center at 5 3
\multiput {$\bullet$} at 1 3  9 3 /
 \put {$P_1$} at .4 3
 \put {$P_2$} at 9.6 3
  \put {$C$} at 7 4.4
 \put {$C'_1$} at 1.7 1
 \put {$C'_2$} at 8.4 1
 \put {$C_j$} at 4.7 .8
 \put {$\Sigma$} at 11.3 5.5
\put {$\longleftarrow$} at 14 3
 \put {$p$} at 14 3.7
\setcoordinatesystem units <.5truecm,.5truecm> point at -19 0
\putrectangle corners at -1.5 -1 and 10.5 6.6
 \putrule from 0 1 to 9 1
 \putrule from 0 3 to 2 3
 \putrule from 0 4.5 to 2 4.5
 \putrule from 1 0 to 1 5.5
 \putrule from 7 3 to 9 3
 \putrule from 7 4.5 to 9 4.5
 \putrule from 8 0 to 8 5.5
 \multiput {$\bullet$} at 1 1 8 1 /
 \put {$P_1$} at 1.6 1.5
 \put {$C'_1$} at -.5 3.3
  \put {$C^{\prime\prime}_1$} at -.5 4.8
  \put {$C$} at -.5 1.3
  \put {$C_1$} at .5 0
  \put {$P_2$} at 7.4 1.5
 \put {$C'_2$} at 9.5 3.3
  \put {$C^{\prime\prime}_2$} at 9.5 4.8
  \put {$C_j$} at 4.2 2.5
  \put {$C_2$} at 8.7 0
  \put {$\tilde \Sigma$} at -2.2 6
\setdashes \setquadratic
  \plot  3.5 2.7  4.5 -.1  5.5 2.7 /
\endpicture
 }
 \vskip .9truecm
  \centerline{\smc Figure 3}

\pagebreak

\noindent {\smc Claim.} {\sl A fibre $C_j$ in $\Sigma$ not
containing $P_1$ or $P_2$ must have $\alpha_j = 1$.} Indeed, by
allowedness, if $\alpha_j \ne 1$ we must have $\alpha_1 = \alpha_2
= 1$ and also $\alpha_\ell = 1$ for possible other fibres $C_\ell$
on $\Sigma$. But then the adjunction formula for $C$ on $\tilde
\Sigma$ yields $\alpha_j = 0$, contradicting the other requirement
for allowedness.

Consequently, for the computation of ${\Cal E}_{\tilde \Sigma}
(p^{-1} D, p^\ast \omega^{1/d})$ we can neglect possible fibres
not containing $P_1$ or $P_2$, and we know that $\alpha_1 +
\alpha_2 = 0$ and $\alpha_1 \ne 0 \ne \alpha_2$. The adjunction
formula on $\tilde \Sigma$ for $C'_i$ and $C^{\prime\prime}_i$
yields $\alpha^\prime_i = \alpha^{\prime \prime}_i =
\frac{\alpha_i + 1}{2}$ for $i = 1,2$. (Here the notations
$\alpha_i,\alpha^\prime_i,\alpha^{\prime \prime}_i$ refer to
Figure 3.) Using the structure of $\Pic \Sigma$ one computes that
the self--intersection $C \cdot C =4$ on $\Sigma$, see e.g. [Ve1,
Remark 6.7]. Then clearly $C \cdot C =0$ on $\tilde \Sigma$.

 So in the defining expression for ${\Cal
E}_{\tilde \Sigma} (p^{-1} D, p^\ast \omega^{1/d})$ we only need
to sum the terms for $I \subset T$ with $\alpha_i \ne 0$ for all
$i \in I$. A simple computation (or directly the formula [Ve5,
Proposition 5.4]) yields
$$\split {\Cal E}_{\tilde \Sigma} (p^{-1} D, p^\ast
\omega^{1/d}) & = L(L-1) + \sum^2_{i=1} \frac{L-1}{L^{\alpha_i}-1}
(L-2+2(1 + L^{\frac{\alpha_i+1}{2}})) \\
& = L(L-1) (1 + \sum^2_{i=1} \frac{1}{L^{\alpha_i}-1}) + 2(L-1)
\sum^2_{i=1} \frac{L^{\frac{\alpha_i+1}{2}}} {L^{\alpha_i}-1} \\
& = 0 + 2(L-1)L^{1/2} \frac{L^{\alpha_1 + \frac{\alpha_2}{2}} -
L^{\frac{\alpha_1}{2}} - L^{\frac{\alpha_2}{2}} + L^{\alpha_2 +
\frac{\alpha_1}{2}}}{(L^{\alpha_1}-1)(L^{\alpha_2}-1)} \\
& = 2(L-1)L^{1/2} \frac{(L^{\frac{\alpha_1 + \alpha_2}{2}} - 1)
(L^{\frac{\alpha_1}{2}} +
L^{\frac{\alpha_2}{2}})}{(L^{\alpha_1}-1)(L^{\alpha_2}-1)} = 0,
\endsplit
$$
using twice that $\alpha_1 + \alpha_2 = 0$. \qed
\enddemo

\bigskip
We are now ready to prove our vanishing theorem.

\bigskip
\proclaim{2.3. Theorem} Let $X$ be a non--singular projective
rational surface and $\omega^{1/d}$ a multi--valued differential
form on $X$ without logarithmic poles. Let $B = \bigcup_{i \in T}
C_i$ be a {\rm connected} normal crossings divisor on $X$
satisfying $B \supset | \div \omega |$. If $\chi(X \setminus B)
\leq 0$, then  \linebreak  $PV \int_X \omega^{1/d}$ $(= {\Cal E}_X
(B, \omega^{1/d})) = 0$.
\endproclaim

\bigskip
\noindent {\it Remark.} The generalization with $B \supset | \div
\omega |$ is not only needed for the applications in \S 3, but is
already useful in the proof for $B = | \div \omega |$.

\bigskip
\demo{Proof} We first explain our strategy. We will construct maps
$\Sigma \overset g \to \leftarrow   \tilde X \overset h \to
\rightarrow X$ as in the Structure Theorem 2.1. If the exceptional
situation of Lemma 1.5(2) does not occur in any blowing--up of $g$
or $h$, then this lemma implies
$${\Cal E}_X(B,\omega^{1/d}) = {\Cal E}_{\tilde X} (h^{-1}B,
\omega^{1/d}) = {\Cal E}_{\tilde X} (h^{-1} B^\prime,
\omega^{1/d}) = {\Cal E}_\Sigma (g(h^{-1}B^\prime),
\omega^{1/d}),$$ where for simplicity we keep the notation
$\omega^{1/d}$ on each surface (a rational differential form is a
birational notion anyway). By Lemma 2.2 this last expression
vanishes (the configuration on $\Sigma$ turns out to be allowed).
The crucial point is that we will show that indeed such an
exceptional situation never occurs.

We will denote all curves appearing in the image of
$h^{-1}B^\prime$ in any intermediate surface by $C_\ell$.

\bigskip
\noindent {\smc Step 1.} {\sl The surface $X$ does not contain any
$(-1)$--curve disjoint from $B$.}

\bigskip
\noindent Indeed, suppose that $A$ is such a $(-1)$--curve
disjoint from $B$. Applying the adjunction formula for $A$ on $X$
would yield $-2 = \deg K_A = A \cdot A = -1$.

\smallskip
So we can really apply Theorem 2.1 and we obtain the extended
configuration $B^\prime \supset B$ on $X$ and the maps $\Sigma
\overset g \to \leftarrow \tilde X \overset h \to \rightarrow X$,
using from now on all notations introduced in that theorem.

\bigskip
\noindent {\smc Step 2.} {\sl The exceptional situation of Lemma
1.5(2) does not occur for {\rm any} constituting blowing--up of $h
: \tilde X \rightarrow X$.}

\bigskip
\noindent Suppose that there is such a blowing--up $b : X^{\prime
\prime} \rightarrow X^\prime$ with centre $P_i \in C_i$, $\alpha_i
= 0$, and two other curves $C_{i_1}$ and $C_{i_2}$ intersecting
$C_i$ outside $P_i$, satisfying $\alpha_{i_1}+\alpha_{i_2}=0$ and
$\{ \alpha_{i_1}, \alpha_{i_2} \} \ne \{-1,1 \}$; see Figure 4.
Denote by $C$ the exceptional curve of $b$; recall that $\alpha =
1$.

\vskip 1truecm
 \centerline{
\beginpicture
\setcoordinatesystem units <.5truecm,.5truecm>
 \putrule from -1 3 to 2.5 3
 \putrule from -1 4.5 to 2.5 4.5
 \putrule from 1 0 to 1 5.5
 \multiput {$\bullet$} at 1 1  1 3  1 4.5 /
 \multiput {$\dots$} at 3.2 3  3.2 4.5 /
 \put {$\vdots$} at 1 6.3
 \put {$P_1$} at 1.7 5
 \put {$P_2$} at 1.7 3.5
 \put {$P_i$} at .3 1
 \put {$C_{i_2}$} at -1.6 3.3
  \put {$C_{i_1}$} at -1.6 4.8
  \put {$C_i$} at 1.5 -.4  \put{$(\alpha_i=0)$} at 3.5 -.3
  \put {$X^{\prime}$} at 6.5 6.4
\put {$\longleftarrow$} at 11 3
 \put {$b$} at 11 3.7
\setcoordinatesystem units <.5truecm,.5truecm> point at -17 0
 \putrule from 0 1 to 5 1
 \putrule from -1 3 to 2.5 3
 \putrule from -1 4.5 to 2.5 4.5
 \putrule from 1 0 to 1 5.5
 \multiput {$\bullet$} at 1 1  1 3  1 4.5 /
 \multiput {$\dots$} at 3.2 3  3.2 4.5 /
 \put {$\vdots$} at 1 6.3
 \put {$P_1$} at 1.7 5
 \put {$P_2$} at 1.7 3.5
 \put {$P_i$} at .4 1.5
 \put {$C_{i_2}$} at -1.6 3.3
  \put {$C_{i_1}$} at -1.6 4.8
  \put {$C$} at  5.6 1.3  \put{$(\alpha=1)$} at 7.4 1.3
  \put {$C_i$} at 1.5 -.4  \put{$(\alpha_i=0)$} at 3.5 -.3
  \put {$X^{\prime\prime}$} at 7.5 6.4
\endpicture
 }
 \vskip 1truecm
  \centerline{\smc Figure 4}
 \vskip 1truecm

Since $b$ is really needed to resolve the indeterminacies of
$\varphi$ ($h$ is minimal), either $C$ itself or some other curve,
created after a chain of blowing--ups starting with centre a point
of $C$, will project by $g$ onto a section or bisection on
$\Sigma$. Such a curve $C^\prime$ has $\alpha^\prime \in \Bbb Z_{>
0}$, and hence the configuration on $\Sigma$ must be case (b).
This also implies that $C_i$ cannot project onto a section or
bisection on $\Sigma$, so it is either blown down during $g$ or
becomes a fibre on $\Sigma$. Both possibilities yield that at some
stage of $g$ the curve $C_i$ can intersect {\sl at most two} other
components of the total transform of $B^\prime$.

So in order to reach that stage by blowing--downs starting from
$\tilde X$, it is necessary that at least for one of the points
$P_1, P_2$ or $P_i$, considered on $C_i \subset \tilde X$, the
following holds (recall that we keep the same notations for strict
transforms of curves and points): {\sl $h^{-1} B^\prime \setminus
\{ P_\ell \}$ has two connected components, and the one not
containing $C_i$ is a tree that gets completely contracted before
reaching that stage of $g$.}

\noindent This cannot happen for $P_i$ because, as explained
above, either $C$ itself or another curve in \lq its tree\rq\ must
project onto a section of $\Sigma$.

Consequently there is such a tree contracting to $P_1$ or $P_2$
during $g$, say to $P_1$. Now, since $\alpha_i = 0$, the component
$C_\ell$ of $h^{-1}B$ through $P_1$ (in $\tilde X$) still has
$\alpha_\ell = \alpha_{i_1}$. Analogously, at the stage of $g$
just before the last blowing--down needed to contract the tree,
the `last' component $C^\prime_\ell$ through $P_1$ still has
$\alpha^\prime_\ell = \alpha_{i_1} (\ne 1)$. But, on the other
hand, since $C^\prime_\ell$ now gets contracted by that
blowing--down, we must have $\alpha^\prime_\ell  = 1$. This
contradiction finishes step 2.

\bigskip
\noindent {\smc Step 3.} {\sl The exceptional situation of Lemma
1.5(2) does not occur for any constituting blowing--up of $g :
\tilde X \rightarrow \Sigma$.}

\bigskip
\noindent (To be precise, in case (c) for $\Sigma$ we consider
only constituting blowing--ups of $\tilde X \to \tilde \Sigma$,
using the notation of (2.2).) Suppose that there is such a
blowing--up $b: X^{\prime \prime} \rightarrow X^\prime$, where we
use again the notations as in Figure 4, but this time `during
$g$'. We consider two cases.

\bigskip (I) {\sl $P_i$, considered on $C_i \subset X^\prime$, does not belong to any other
component of the total transform of $B^\prime$.}

\medskip
 \noindent Here we consider two subcases.

\medskip
(i) $C \subset h^{-1}B$ in $\tilde X$.

\medskip
 \noindent Since $C_i(\subset \tilde X)$ must be blown down
during $h$ (all $\alpha_\ell \ne 0$ on $X$) and $B$ has normal
crossings, {\sl at most two} components of the total transform of
$B$ can intersect $C_i$ just before that blowing--down. As we
explained above, since $\alpha_{i_1} \ne 1 \ne \alpha_{i_2}$, the
components intersecting $C_i$ in $P_1$ and $P_2$ cannot be blown
down before $C_i$. So the tree that was created during $g$,
starting with the blowing--up with centre $P_i$, must be
contracted again by $h$. This contradicts the minimality of $h$.

\medskip (ii) $C(\subset \tilde X)$ is a component of (the strict
transform of) $B^\prime \setminus B$.

\medskip
 \noindent Because $h^{-1}B$ is connected, each possible
further blowing--up of $g$ (after $b$) with centre in $C$ must
have $P_i$ as centre. Consequently
$$\chi(X \setminus B^\prime) < \chi (X \setminus B) \leq 0$$
and $\Sigma$ must be as in case (a) with {\sl at least three
fibres}. The strict transform in $\tilde X$ of the section $C_1
\subset \Sigma$ can intersect {\sl at most two} other components
of $h^{-1}B$. (Indeed, since $\alpha_i=0$, the morphism $h$ cannot
be the identity. Hence $C_1$, being the only non--fibre in
$\Sigma$, must be created by the last blowing--up of $h$.)
 So $C_1$ must intersect in
$\tilde X$ at least one component $C^\prime$ coming from $B^\prime
\setminus B$. But $C^\prime$ gets contracted during $g$; so in
order to have at least three intersections with $C_1$ in $\Sigma$,
$C^\prime$ must intersect in $\tilde X$ {\sl another} component of
$h^{-1}B$. This contradicts Theorem 2.1(3).

\bigskip
(II) {\sl $P_i$, considered on $C_i \subset X^\prime$, also
belongs to another component $C^\prime$ of the total transform of
$B^\prime$. Then this curve $C'$ should be drawn (in the reader's
mind) through $P_i$ on the $X'$--side of Figure 4, and its strict
transform through another point of $C$ on the
$X^{\prime\prime}$--side.}

\medskip
 \noindent Since, by Lemma 1.5, $(h^{-1}B, \omega^{1/d})$
is allowed on $\tilde X$, we then have necessarily $\alpha^\prime
= 1$. Moreover, we claim that $C_i$ must already live on $\Sigma$,
or on $\tilde \Sigma$ in case (c), where we use the notation of
(2.2). Indeed, suppose that $C_i$ is created by a blowing--up of
$g$, or of $\tilde X \to \tilde \Sigma$ in case (c); at that stage
$C_i$ can intersect {\sl at most two} other components. Thus $C'$
should be created afterwards during $g$. But this situation is
precisely the previous case (I)(ii), which we already
contradicted.

We now investigate the cases (a), (b) and (c) for $\Sigma$ in our
(hypothetical) situation. Case (a) is clearly impossible. In case
(b) $C_i$ must be a section, and in case (c) $C_i$ must be the
bisection in $\Sigma$, see Figure 5. (In case (c) $C_i$ could be a
priori another curve in $\tilde \Sigma$, but then we would have
two intersecting curves with $\alpha_\ell = 0$ on $\tilde \Sigma$,
implying the same phenomenon on $\tilde X$.) Somewhat abusing
notation in Figure 5, the curves $C_{i_1}$, $C_{i_2}$ and $C'$ in
$\Sigma$ could be different from those intersecting $C_i$ in $X'$,
but their $\alpha$--coefficients are also $\alpha_{i_1}$,
$\alpha_{i_2}$ and $\alpha' (=1)$, respectively.

\vskip 1truecm
 \centerline{
\beginpicture
\setcoordinatesystem units <.5truecm,.5truecm>
 \putrectangle
corners at -.5 0 and 13 6
 \putrule from 2 1 to 2 5
 \putrule from 6.5 1 to 6.5 5
  \putrule from 4 1 to 4 5
  \putrule from 1 2 to 8 2
 \putrule from 1 4 to 8 4
 \multiput {$\bullet$} at 2 2  4 2  6.5 2  2 4  4 4  6.5 4   /
 \put {$P_1$} at 1.5 4.6
 \put {$P_i$} at 6 4.6
 \put {$P_2$} at 3.5 4.6
 \put {$P'_1$} at 1.5 2.6
 \put {$P'_i$} at 6 2.6
 \put {$P'_2$} at 3.5 2.6
  \put {$C_i$} at 8.7 4
  \put {$C'_i$} at 8.7 2
  \put {$(\alpha_i=0)$} at 10.7 4
  \put {$(\alpha'_i=0)$} at 10.7 2
 \put {$C_{i_1}$} at 1.7 .7
 \put {$C_{i_2}$} at 3.7 .7
 \put {$C'$} at 6 .7
 \put {$\Sigma$} at 13.8 5.5
 \put {$(b)$} at -2 3
\setcoordinatesystem units <.5truecm,.5truecm> point at 0 9
\putrectangle corners at -.5 0 and 10.5 6
 \putrule from 1 1 to 1 5
 \putrule from 9 1 to 9 5
  \putrule from 4 1 to 4 5
\ellipticalarc axes ratio 4:1  360 degrees from 1 3 center at 5 3
\multiput {$\bullet$} at 1 3  9 3  4 3.95 /
 \put {$P_1$} at .4 3
 \put {$P_2$} at 9.7 3
 \put {$P_i$} at 3.5 4.5
  \put {$C_i$} at 6.7 3.4
  \put {$(\alpha_i=0)$} at 6.7 4.5
 \put {$C'$} at 4.6 .8
 \put {$\Sigma$} at 11.2 5.5
\put {$\longleftarrow$} at 13.5 3
 \put {$p$} at 13.5 3.7
 \put {$(c)$} at -2 3
\setcoordinatesystem units <.5truecm,.5truecm> point at -18.5 9
\putrectangle corners at -2.5 -1 and 10.5 6.6
 \putrule from -1 1 to 9 1
 \putrule from 0 3 to 2 3
 \putrule from 0 4.5 to 2 4.5
 \putrule from 1 0 to 1 5.5
 \putrule from 7 3 to 9 3
 \putrule from 7 4.5 to 9 4.5
 \putrule from 8 0 to 8 5.5
 \multiput {$\bullet$} at 1 1  8 1  3.9 1 /
 \put {$P_1$} at 1.6 1.5
 \put {$P_i$} at 3.4 .4
  \put {$C_i$} at -.6 1.5
  \put {$(\alpha_i=0)$} at -.7 .4
  \put {$C_{i_1}$} at 1.7 5.7
  \put {$P_2$} at 7.5 1.6
  \put {$C'$} at 4.9 2.5
  \put {$C_{i_2}$} at 8.7 5.7
  \put {$\tilde \Sigma$} at -3.2 6
 \setquadratic
  \plot  3.5 2.6  4.5 -.1  5.5 2.6 /
\endpicture
 }
 \vskip 1truecm
  \centerline{\smc Figure 5}
 \vskip 1truecm

\noindent Case (b). On $\tilde X$ we have that $C_i$ (being an
exceptional curve of $h$) can intersect {\sl at most two} other
components of $h^{-1}B$. So all but at most two intersections of
$C_i$ in $\Sigma$ with fibres $C_\ell$ of $\pi : \Sigma
\rightarrow \Bbb P^1$ must `split' during $g$, i.e. $C_i$ must
eventually get separated from those $C_\ell$ by a blowing--up of
$g$ with centre in $C_i$ and exceptional curve coming from
$B^\prime \setminus B$. (The same is true for $C'_i$.)

 In order to obtain such an exceptional curve
$C_j$ from $B^\prime \setminus B$, intersecting $C_i$ in $P_1$, we
should have simultaneously $\alpha_j = 1$ and $\alpha_j =
\alpha_{i_1} (\ne 1)$, by the same argument as for the conclusion
of step 2. The same contradiction works for $P_2, P^\prime_1$ and
$P^\prime_2$. Hence we must separate $C_i$ and $C^\prime_i$ from
$C^\prime$ with exceptional curves from $B^\prime \setminus B$,
intersecting $C_i$ and $C^\prime_i$ in $P_i$ and $P^\prime_i$,
respectively. But this contradicts the connectivity of $h^{-1}B$.

\smallskip \noindent Case (c). Analogously we will have to separate
$C^\prime$ from $C_i$, contradicting the connectivity of
$h^{-1}B$.

\bigskip
\noindent {\smc Step 4.} {\sl $(g(h^{-1}B^\prime),\omega^{1/d})$
is allowed on $\Sigma$ and ${\Cal E}_X (B,\omega^{1/d}) = {\Cal
E}_\Sigma(g(h^{-1}B^\prime),\omega^{1/d}) = 0.$}

\bigskip
\noindent
By Lemma 1.5 we have ${\Cal E}_X(B,\omega^{1/d}) = {\Cal
E}_{\tilde X}(h^{-1}B, \omega^{1/d})$. Since also $h^{-1}B^\prime$
is a normal crossings divisor and the curves $C_\ell$ in $B^\prime
\setminus B$ have $\alpha_\ell = 1$, we have that also
$(h^{-1}B^\prime, \omega^{1/d})$ is allowed on $\tilde X$ and
${\Cal E}_{\tilde X}(h^{-1}B, \omega^{1/d}) = {\Cal E}_{\tilde
X}(h^{-1}B^\prime, \omega^{1/d})$. It is easy to see that then
necessarily $(g(h^{-1}B^\prime),\omega^{1/d})$ is allowed on
$\Sigma$, and thus ${\Cal E}_{\tilde X}(h^{-1}B^\prime,
\omega^{1/d})={\Cal E}_\Sigma(g(h^{-1}B^\prime),\omega^{1/d})$,
again by Lemma 1.5(2). This last expression equals zero by Lemma
2.2.
 \qed
\enddemo

\bigskip
\noindent {\bf 2.4.} {\sl Remark.} Theorem 2.3 and its proof are
quite subtle. In [Ve7, 3.4] we constructed the following similar
example. Let $X = \Bbb P^2$ and $\omega^{1/2}$ a multi--valued
differential form with $|\div \omega|$ a (non--singular) conic
$B$. So $\omega^{1/2}$ has no logarithmic poles on  $X$ and one
easily computes that $PV \int_X \omega^{1/2} \ne 0$. Taking
$B^\prime$ as the union of $B$ and one tangent line to $B$, we can
construct $\Sigma \overset g \to \leftarrow \tilde X \overset h
\to \rightarrow X$, where $h$ is a composition of three
blowing--ups and $g$ a composition of two blowing--downs, all
outside of $X \setminus B^\prime$, such that $g(h^{-1}B^\prime)
\subset \Sigma$ is case (b) in Theorem 2.1 with exactly one fibre.
Also on $\Sigma$ we have that $\omega^{1/2}$ has no logarithmic
poles, but now $PV \int_\Sigma \omega^{1/2} = 0$. In this example
we have that $\chi(X \setminus B^\prime) = 0$; the only
obstruction with the data of Theorem 2.3 is that here $\chi(X
\setminus B) = 1 \, (> 0)$. And in fact the `change' in principal
value integral is caused by $g$ which consists precisely of the
exceptional situation of Lemma 1.5(2).

Another example in this context is provided by the minimal
embedded resolution of two smooth plane conics with one
intersection point, see [ACLM1, Example 2.14].

\bigskip
\noindent {\bf  2.5.} Our vanishing theorem specializes to the
level of Hodge polynomials or Euler characteristics. With the same
notations as in Definition 1.4 we can introduce analogously the
invariants
$$E_X(D,\omega^{1/d}) := \sum_{\Sb I \subset T \\ \forall i \in
I:\alpha_i \ne 0 \endSb}  H(C^\circ_I) \prod_{i \in I}
\frac{uv-1}{(uv)^{\alpha_i}-1} + \sum_{\Sb i \in T \\ \alpha_i=0
\endSb} (-C_i \cdot C_i) \prod_{j \in T_i} \frac{uv - 1}{(uv)^{\alpha_j}-1}$$
and
$$e_X(D,\omega^{1/d}) := \sum_{\Sb I \subset T \\ \forall i \in I
: \alpha_i \ne 0 \endSb} \chi (C^\circ_I) \prod_{i \in I}
\frac{1}{\alpha_i} + \sum_{\Sb i \in T \\ \alpha_i = 0 \endSb}
(-C_i \cdot C_i) \prod_{j \in T_i} \frac{1}{\alpha_j}.$$ Then with
the same data as in Theorem 2.3 we obtain that
$E_X(B,\omega^{1/d}) = e_X (B,\omega^{1/d}) = 0$.

\bigskip
\noindent {\bf 2.6.} Let $L$ be a $p$--adic field with valuation
ring ${\Cal O}$, maximal ideal $P$ and residue field $\frac{\Cal
O}{P} \cong \Bbb F_q$. We choose an embedding of $L$ into $\Bbb
C$. When $X$ and $\omega^{1/d}$ are defined over $L$, we can
introduce analogously
$$E^L_X(D,\omega^{1/d}) := \sum_{\Sb I \subset T \\ \forall i \in
I:\alpha_i \ne 0 \endSb} \card (C^\circ_I)_{\Bbb F_q} \prod_{i \in
I} \frac{q-1}{q^{\alpha_i}-1} + \sum_{\Sb i \in T \\ \alpha_i = 0
\endSb} (-C_i \cdot C_i) \prod_{j \in T_i}
\frac{q-1}{q^{\alpha_j}-1},$$ where $\card(\cdot)_{\Bbb F_q}$
denotes the number of $\Bbb F_q$--rational points of the reduction
mod $P$ of $C^\circ_I$. When $\omega^{1/d}$ has no logarithmic
poles and if suitable conditions about good reduction mod $P$ are
satisfied, a similar proof as for Denef's formula for the
$p$--adic Igusa zeta function [De1] yields that
$E^L_X(D,\omega^{1/d})$ is (up to a power of $q$) precisely the
$p$--adic principal value integral $PV \int_{X(L)} | \omega^{1/d}
|$.

When we assume that the field $L$ is `big enough', the same data
as in Theorem 2.3 will again imply that this invariant vanishes.
More precisely we want that all constructions in this proof of the
theorem can be done `over $L$', for instance that all centres of
blowing--ups are $L$--rational points. This can always be achieved
by taking a finite extension of a given $p$--adic field.

\bigskip
\bigskip
\noindent {\bf 3. Cancelation of candidate poles for zeta
functions}

\bigskip
\noindent {\bf 3.1.} Fix a polynomial $f \in \Bbb
C[X_1,\dots,X_{n+1}] \setminus \Bbb C$, or, more generally, a
non--constant regular function $f : M \rightarrow \Bbb A^1$ from a
non--singular $(n+1)$--dimensional variety $M$. Denef and Loeser
[DL2] associated to $f$ its {\sl motivic zeta function} ${\Cal
Z}_{\mot} (f;T) \in {\Cal M}[[T]]$, where ${\Cal M}$ is the
localization of $K_0(Var)$ with respect to $L$. It describes the
orders of all truncated arcs on $M$ along the hypersurface $\{ f =
0 \}$, and is in fact modeled on the classical $p$--adic Igusa
zeta function. (See also [DL4] or [Ve6] for an introduction to
this topic.) We just mention a formula for $Z_{\mot}(f;T)$ in
terms of an embedded resolution of $\{ f = 0 \}$, implying in
particular that this invariant belongs to the localization of
${\Cal M}[T]$ with respect to $L^a - T^b$, where $a,b \in \Bbb
Z_{> 0}$.

\bigskip
\noindent \proclaim {3.2. Theorem {\rm [DL2, 2.2.1]}} Let $h : Y
\rightarrow M$ be an embedded resolution of $\{ f = 0 \}$. Let
$E_j, j \in K$, be the irreducible components of $h^{-1} \{ f = 0
\}$. Denote by $N_j$ the multiplicity of $E_j$ in $\div(f \circ
h)$ and by $\nu_j - 1$ the multiplicity of $E_j$ in $\div(h^\ast
dx)$, where $dx$ is a local generator of the sheaf of differential
$(n+1)$--forms on $M$. For $J \subset K$ we put $E^\circ_J :=
(\bigcap_{j \in J} E_j) \setminus (\bigcup_{\ell \not\in J}
E_\ell)$; so $Y = \coprod_{J \subset K} E^\circ_J$. Then
$${\Cal Z}_{\mot} (f;T) = L^{-(n+1)} \sum_{J \subset K} [E^\circ_J]
\prod_{j \in J} \frac{(L-1)T^{N_j}}{L^{\nu_j}- T^{N_j}}.$$
\endproclaim

\bigskip
\noindent {\bf 3.3.} This zeta function specializes to the {\sl
Hodge zeta function} $Z_{\Hod}(f;T)$, replacing in the formula
above all classes of varieties in $\Cal M$ by their Hodge
polynomial, and further to the `classical' {\sl topological zeta
function}
$$z_{\top}(f;s) = \sum_{J \subset K} \chi(E^\circ_J) \prod_{j \in
J} \frac{1}{\nu_j + sN_j} \in \Bbb Q(s)$$ of [DL1]. We refer to
[DL2] or [Ve6] for more details.

\bigskip
\noindent {\bf 3.4.} The famous {\sl monodromy conjecture}, stated
originally for the $p$--adic Igusa zeta function, can be
formulated for ${\Cal Z}_{\mot}(f;T)$ as follows [DL2, 2.4]:

\bigskip
\noindent {\sl $Z_{\mot}(f;T)$ belongs to the localization of
${\Cal M}$ with respect to those $L^a - T^b$,  $a,b \in \Bbb Z_{>
0}$, such that $e^{2 \pi ia/b}$ is an eigenvalue of the local
monodromy on the Milnor fibre of $f$ at some point of $\{ f = 0
\}$.}

\bigskip
\noindent So if $L^{\nu/N}$ is a pole of ${\Cal Z}_{\mot}(f;T)$,
then $e^{2 \pi i \nu/N}$ is expected to be an eigenvalue of the
local monodromy. Note however that one has to be careful with the
notion of pole here, the difficulty being that we do not know
whether ${\Cal M}$ is a domain. See e.g. [RV2] for a precise
definition. For the Hodge and topological zeta function the notion
of pole is clear.

The ($p$--adic version of) the conjecture was proved in dimension
two $(n=1)$ by Loeser [Loe1]; a simple proof in the
motivic/Hodge/topological setting is in [Ro1]. It is still open in
general, with partial results in dimension three [Ve2,RV1,ACLM1]
in and other special cases [Loe2,ACLM2].

\bigskip
\noindent {\bf 3.5.} We keep using the notation of Theorem 3.2.
Fix an exceptional component $E_j$ of $h$ which is mapped to a
point by $h$. It induces the candidate pole $L^{\nu_j/N_j}$ for
${\Cal Z}_{\mot}(f;T)$, respectively $(uv)^{\nu_j/N_j}$ for
$Z_{\Hod}(f;T)$ and $-\nu_j/N_j$ for $z_{\top}(f;s)$.

In order for the monodromy conjecture to hold, looking at
A'Campo's formula for the monodromy zeta function [A'C] one
expects the following. Suppose we are in the generic case that
$\nu_j/N_j \ne \nu_i/N_i$ for all $i \ne j$. If $\chi(E^\circ_j) =
0$, maybe even if $(-1)^n \chi (E_j^\circ) \leq 0$, then `in
general' $L^{\nu_j/N_j}$ should not be a pole of ${\Cal
Z}_{\mot}(f;T)$.

Somewhat more precise, suppose only that $\nu_j/N_j \ne \nu_i/N_i$
for all $i \in S_j := \{ i \in K \mid E_i \cap E_j \ne \emptyset
\}$. Then, if $(-1)^n \chi (E^\circ_j) \leq 0$, one expects that
`in general' $E_j$ does not contribute to the possible pole
$L^{\nu_j/N_j}$, which means that $L^{\nu_j/N_j}$ should not be a
pole of

$$\frac{1}{L^{n+1}} \sum_{j \in I \subset K} [E^\circ_I] \prod_{i
\in I} \frac{(L-1)T^{N_i}}{L^{\nu_i}-T^{N_i}}. \leqno (\ast)$$

\noindent We refer to e.g. [Ve2, \S1] for a motivation. For $n=1$
this expectation is true and is part of the proof of the monodromy
conjecture [Loe1,Ro1]. Now, $L^{\nu_j/N_j}$ not being a pole of
($\ast$) can be reformulated as `its residue is zero', i.e.
$$\sum_{j \in I \subset K} [E^\circ_I] \prod_{i \in I \setminus \{
j \}} \frac{L-1}{L^{\alpha_i}-1} \leqno (\ast \ast)$$ is zero,
where $\alpha_i := \nu_i - (\nu_j/N_j) N_i$ for $i \in S_j$. Now
($\ast \ast$) is a motivic principal value integral on $E_j$.
Indeed, let $dx$ be a local generator of the sheaf of
$(n+1)$--forms on $M$ around the point $h(E_j)$. Then the
Poincar\'e residue $\omega^{1/d}$ of $(f \circ h)^{-\nu_j/N_j}
h^\ast(dx)$ on $E_j$ is a multi--valued differential form on $E_j$
with $\div \omega^{1/d} = \sum_{i \in S_j} (\alpha_i - 1) (E_j
\cap E_i)$. This is easily verified with local coordinates, see
also [Ja3]. We thus have that ($\ast \ast$) is (up to a non--zero
constant) equal to $PV \int_{E_j} \omega^{1/d}$.

Note however that by assumption all $\alpha_i \ne 0$, so
$\omega^{1/d}$ indeed has no logarithmic poles, but it is possible
that some $\alpha_i = 1$. So in fact we land in a natural way in
the more general framework of (1.3) with $D := \bigcup_{i \in S_j}
(E_j \cap E_i) \supset | \div \omega |$, where the inclusion may
be strict, and $(\ast \ast)$ is (essentially) ${\Cal E}_{E_j} (D,
\omega^{1/d})$.

\bigskip
\noindent {\bf 3.6.} Let now $n=2$. We may assume that $h$ is
constructed as a composition of blowing--ups with non--singular
centre. Fix as above a projective exceptional surface $E_j$ which
is mapped to a point by $h$. The surface $E_j$ was created during
some blowing--up $\pi$ of the resolution process $h$ as a surface
$E^{(0)}_j$, where either $E^{(0)}_j \cong \Bbb P^2$ or
$E^{(0)}_j$ is a ruled surface, when the centre of $\pi$ is a
point or a curve, respectively. And then $E_j$ is obtained from
$E^{(0)}_j$ by a composition $\varphi : E_j \rightarrow E^{(0)}_j$
of (point) blowing--ups. Denote again $D := \bigcup_{i \in S_j}
(E_j \cap E_i)$. We have that $D$ is the inverse image by
$\varphi$ of the intersection of $E^{(0)}_j$ with the other
components of the total inverse image of $\{ f = 0 \}$ at the
stage of $h$ when $E^{(0)}_j$ was just created. In particular $D$
is connected if and only if this intersection on $E^{(0)}_j$ is
connected. And it is thus always connected if $E^{(0)}_j \cong
\Bbb P^2$.

\bigskip
\noindent {\bf 3.7.} By the considerations above Theorem 2.3
yields the following cancelation result for candidate poles of the
motivic zeta function. For ease of reference we recall the
notations.

Let $M$ be a three--dimensional non--singular variety and $f : M
\rightarrow \Bbb A^1$ a non--constant regular function. Let $h : Y
\rightarrow M$ be an embedded resolution of $\{ f = 0 \}$,
constructed as a composition of blowing--ups. Denote by $E_j,j \in
K$, the irreducible components of $h^{-1} \{ f = 0 \}$ and let
$N_j,  \nu_j$ and $E^\circ_J$ be as in (3.2). Suppose that $E_j$
is mapped to a point by $h$ and that $\nu_j/N_j \ne \nu_i/N_i$ for
all $i \in S_j := \{ i \in K \mid E_j \cap E_i \ne \emptyset \}$.
Denote
$$R_{E_j} := \sum_{j \in I \subset K} [E^\circ_I] \prod_{i \in I
\setminus \{ j \}} \frac{L-1}{L^{\alpha_i}-1} ,$$
 `the contribution of $E_j$ to the residue of $L^{\nu_j/N_j}$ for ${\Cal
Z}_{\mot}(f;T)$'.

\bigskip
\proclaim{Theorem} Let $\chi(E_j^\circ)\leq 0$.

(1) If $E_j$ is created by blowing up a point, then we have {\rm
always} $R_{E_j} = 0$.

(2) If $E_j$ is created by blowing up a rational curve, and if
$\bigcup_{i \in S_j}(E_j \cap E_i)$ is connected, then $R_{E_j}  =
0$.
\endproclaim

\bigskip
\noindent Recall that in case (2) this connectivity is equivalent
to the connectivity of the analogous intersection configuration
$D^{(0)}$ on the (rational) ruled surface $E^{(0)}_j$ (3.6). Note
then that the exceptions in (2), i.e. those $E_j$ with a
non--connected intersection configuration, are very special!
Indeed, non--connectivity of $D^{(0)}$ implies for instance that
$D^{(0)}$ does not contain any fibre of the ruled surface
$E_j^{(0)}$. In the embedded resolution process this is quite
rare.

There is a recent vanishing result of Rodrigues [Ro2] in this
exceptional case. When $\chi(E_j^\circ) = 0$, and assuming a minor
extra condition, he classified all possible non--connected
$D^{(0)}$ with non--singular irreducible components, and verified
that then again $R_{E_j} = 0$.

\bigskip
\noindent {\bf 3.8.} The case where $E_j$ is a rational surface is
the most difficult one. There is no classification of the possible
intersection configurations on $E_j$ with $\chi(E^\circ_j) \leq
0$; instead we used our structure theorem 2.1. When $E_j$ is
created by blowing up a non--rational curve, we already obtained a
classification of the possible configurations with
$\chi(E^\circ_j) \leq 0$ in [Ve2] and [Ve3].

\bigskip
\proclaim{Proposition} We use all notations of (3.7). Let
$\chi(E^\circ_j) \leq 0$. If $E_j$ is created by blowing up a
non--rational curve, then we have always $R_{E_j} = 0$.
\endproclaim

\medskip
 \demo{Proof} Let as above $E_j$ be created as the ruled
surface $E^{(0)}_j$ while blowing up a curve of genus $g$ during
$h$. When $g \geq 2$, we classified the possible configurations on
$E^{(0)}_j$ with $\chi(E^\circ_j) \leq 0$ in [Ve2, Proposition
5.13], and verified in [Ve2, Propositions 5.1 and 5.3] that
$R_{E_j} = 0$ for them. (The calculation there is in the context
of $p$--adic Igusa zeta functions, but is essentially the same in
the motivic setting.)

When $g = 1$, the possible configurations ${\Cal C}$ of curves on
the ruled surface $E^{(0)}_j$ with $\chi(E_j \setminus {\Cal C})
\leq 0$ were classified in [Ve3, Theorem 6.5]. Again by [Ve2,
Propositions 5.1 and 5.3] we have that $R_{E_j} = 0$, except for
one annoying case (where $R_{E_j} \ne 0$). More precisely, in this
case the curves on the ruled surface consist of a number of
disjoint elliptic curves, where either one of the curves is not a
section, or there are at least three curves. Now recently
Rodrigues showed in [Ro2] that in fact this configuration cannot
occur in the context of exceptional surfaces in an embedded
resolution. \qed
\enddemo

\bigskip
\noindent {\bf 3.9.} The theorem and proposition above provide a
strong confirmation for the Monodromy Conjecture for surfaces
($n=2$), and could be a major contribution to a proof. Of course
there are still various non--obvious remaining parts within this
strategy, for instance extending the results of [Ro2] to {\sl all}
possible non--connected $D^{(0)}$ and handling candidate poles of
higher order.

\bigskip
\noindent {\bf 3.10.} Theorem 3.7 and Proposition 3.8 specialize
to the analogous results in the context of the Hodge and
topological zeta functions of (3.3). They are also valid in the
context of $p$--adic Igusa zeta functions (see e.g. [De1,Ve2]) if
the $p$--adic field $L$ is assumed `big enough' as in (2.6), and
if suitable conditions concerning good reduction mod $P$ as in
[De1] are satisfied. Alternatively, one can take a big enough
number field $F$. Then our vanishing results will be true for the
Igusa zeta functions over all except a finite number of
completions $L$ of $F$.


\enddocument